\documentclass[11pt]{amsart}

\usepackage[utf8]{inputenc}
\usepackage{fullpage,url,amssymb,enumitem,colonequals, mathrsfs, comment}
\usepackage{microtype}
\usepackage{hyperref}
\usepackage{xcolor}
\usepackage{lineno}
\usepackage{cleveref}

\usepackage[pdftex]{}

\newcommand{\C}{\mathbb{C}}
\newcommand{\Cb}{\mathbb{C}}

\newcommand{\F}{\mathbb{F}}

\newcommand{\N}{\mathbb{N}}
\newcommand{\Nb}{\mathbb{N}}

\newcommand{\Q}{\mathbb{Q}}
\newcommand{\Qb}{\mathbb{Q}}
\newcommand{\R}{\mathbb{R}}
\newcommand{\Z}{\mathbb{Z}}
\newcommand{\Zb}{\mathbb{Z}}
\newcommand{\Qbar}{{\overline{\Q}}}

\newcommand{\Ol}{\mathcal{O}}

\newcommand{\eps}{\varepsilon}

\newtheorem{theorem}{Theorem}[section]
\newtheorem{lemma}[theorem]{Lemma}
\newtheorem{corollary}[theorem]{Corollary}
\newtheorem{proposition}[theorem]{Proposition}

\theoremstyle{remark}
\newtheorem{remark}[theorem]{Remark}

\theoremstyle{definition}

\newtheorem{conjecture}[theorem]{Conjecture}
\newtheorem{example}[theorem]{Example}

\usepackage{microtype}

\title{Root sharing and $2$-adic interpolation for Legendre polynomials}
\author{Borys Kadets}
\address{Einstein Institute of Mathematics, Edmond J. Safra campus, Hebrew University, Jerusalem 9190401, Israel}
\email{kadets.math@gmail.com}
\urladdr{\url{http://bkadets.github.io}}

\author{Dan Mangoubi}
\address{Einstein Institute of Mathematics, Edmond J. Safra campus, Hebrew University, Jerusalem 9190401, Israel}
\email{dan.mangoubi@mail.huji.ac.il}
\urladdr{\url{http://www.ma.huji.ac.il/~mangoubi/}}

\begin{document} 
\begin{abstract}
Stieltjes's conjecture on Legendre polynomials asserts that no two distinct Legendre polynomials can share a nonzero root. This question is, surprisingly, connected to a number of disparate mathematical areas. From the point of view of spectral geometry, this question concerns zero-sets of spherical harmonics. A disk counterpart to the conjecture is a famous theorem of Siegel on nonexistence of nonzero common roots of Bessel functions, proved using transcendence theory. A weak version of the conjecture is implied by the dynamical Mordell--Lang conjecture. We prove, using a $2$-adic analytic approach \`a la Skolem,  that for a root $\alpha\neq 0$ of $P_l$ the set $\{n\in \N \colon P_n(\alpha)=0\}$ has cardinality at most $2l$. In particular, this set is finite.
\end{abstract}

\maketitle
\section{Introduction}
In his 1890 letter to Hermite \cite{stieltjes-hermite-letter}, Stieltjes raised the following conjecture.

\begin{conjecture}[Stieltjes]\label{conj:stieltjes}
    Two distinct Legendre polynomials $P_n$ cannot share a nonzero root.
\end{conjecture}

In this paper we prove, using a $2$-adic analytic method, the following result.

\begin{theorem}\label{main-theorem-intro}
    Suppose $\alpha\neq 0$ is a root of the Legendre polynomial $P_l$. Then the set 
    \[\{n\in \N \colon P_n(\alpha)=0\}\] has cardinality at most $2l$.
\end{theorem}
We remark that even the finiteness of the set considered in~Theorem~\ref{main-theorem-intro} has not been previously known.
\Cref{conj:stieltjes}, while appearing very special, arises in a number of contexts in mathematical analysis, number theory, and arithmetic dynamics, which we now discuss.

\subsection{Spectral geometry}
Legendre polynomials are closely connected to ($2$-dimensional) spherical harmonics. Stieltjes's conjecture would show that the nodal sets of zonal spherical harmonics aligned with respect to the axis of symmetry do not have nonempty intersections except maybe at the equator. This formulation connects the conjecture to a number of other results and problems. 

The celebrated analogous problem from the 1860s, concerning common roots of Bessel functions, is known as Bourget's hypothesis and was solved by Siegel~\cite{siegel1929, siegel-book, siegel1929-translation} in 1929 using transcendental number theory. It shows that no two Bessel functions can share a nonzero root. In particular it follows that Dirichlet eigenvalues of the planar disk have no nontrivial multiplicities and that the nodal domains pattern of any eigenfunction on the disk looks like a circular chess-board.

A similar question for the standard $2d$ flat torus has been studied by Bourgain and Rudnick~\cite{bourg-rud2011nodal}. They show that only a finite number of eigenfunctions can vanish on a given curve of nonzero geodesic curvature. 

An instance of the analogous spherical question was studied in \cite{mangoubi-weller-weiser2019harmonic}, which connects also to rigidity questions for  sets of harmonic functions defined by their zeros. The $2$-adic methods used in that paper led us to the present work.

\subsection{Rubel's conjecture}

Rubel \cite[Problem 16]{rubel1992problems} raised the following interesting conjecture. Suppose $f(x)$ is a solution of a linear differential equation with polynomial coefficients, and suppose~$f$ has a Taylor expansion at~$0$: $f(x)=\sum_k a_k x^k$. Then the set~$\{n\in \N\colon a_n=0\}$ is a finite union of arithmetic progressions. This conjecture is a vast generalization of the classical Skolem--Mahler--Lech theorem (see e.g. \cite[pp. 312-318]{zannier2024} and~\cite[ch. 4]{borevich-shafarevich66}), which handles the case of rational functions $f(x)$. The connection to Stieltjes's conjecture comes from the formula

\[\left(1-2xt+x^2\right)^{-1/2}=\sum_{l=0}^\infty P_l(t)x^l.\]
Indeed, $f(x)=(1-2xt+x^2)^{-1/2}$ satisfies
\[(1-2xt+x^2)f'(x)+(x-t)f(x)=0.\]
Thus our \Cref{main-theorem-intro} confirms, in this special case, Rubel's conjecture. There is a vast literature on studying Rubel's conjecture, often formulated in the equivalent language of polynomial--linear recurrences; see, for example, \cite{methfessel2000zeros, bell2006gen-sml, bell2008corrigendum, bell2012zero, wibmer2015sml}. This question is closely related to the wide open dynamical Mordell--Lang conjecture; see the discussion of the connection and of the partial results on Rubel's conjecture in \cite[\S 3.2.1]{bell2016dynamical}.

\subsection{Arithmetic of classical polynomials}

Algebraic number theory provides powerful tools for studying irreducibility and Galois groups of polynomials. For a vast majority of classical polynomial families, one can often combine local information, computed using the theory of Newton polygons, to conclude that polynomials in a given family are irreducible and have Galois group $S_n$ (or occasionally $A_n$, or $S_n\wr S_2$, when the polynomials in the family are even or odd); see \cite{coleman1987galois} for an introduction to this method. This is the case for Hermite polynomials \cite{schur1931affektlose}, for generalized Laguerre polynomials \cite{hajir2005galois}, for Bessel polynomials \cite{filaseta1995irreducibility}, \cite[Chapter 12]{grosswald2006bessel}. Among such classical families the Legendre family is the odd one out, and its arithmetic is still far from being understood, as evidenced by \Cref{conj:stieltjes}. For more on Galois groups of Legendre polynomials see \cite{cullinan2014galois}.

\subsection{Description of our approach}
The main discovery of this paper is an analyticity property of Legendre Polynomials as functions of their~\emph{degree}. Namely, we  consider $P_l(\alpha)$ as a function of~$l$ rather than~$\alpha$. It turns out that after a suitable nonlinear change of variables introduced by  Holt~\cite{holt1913legendre} $P_l(\alpha)$ becomes $2$-adic analytic in~$l$. This idea is similar to the classical method introduced by Skolem in his theorem  on the zeroes of linear recurrent sequences \cite{skolem1934}.
Once we have analyticity we can apply general finiteness results analogous to ones of classical complex analysis results, e.g., a nonzero analytic function has only a finite number of zeros in a closed disk.
 We believe that this analyticity property will find generalizations to other special functions arising in representation theory.
 In more detail,
 Holt introduces the polynomials $H_l(z)$: 
\[H_l(z)\colonequals \binom{2l}{l}^{-1}(z-4)^lP_{2l}\left(\left(\frac{z}{z-4}\right)^{1/2}\right).\]

A similar transformation can be made for odd degree Legendre polynomials; we will restrict ourselves to the even degree case in this discussion.  Our main result is then the following:
\begin{theorem}\label{interpolation-intro}
    Suppose $\beta \in \overline{\Q}_2$ is a root of $H_n(z)$ for some fixed $n$. Consider $\check{H}_{\beta}(l):=H_{l}(\beta)$ as a function from $\N$ to $\Q_2(\beta)$. Then $\check{H}_{\beta}$ extends to a locally-analytic $\Q_2(\beta)$-valued function $\check{H}_{\beta}$ on $\Z_2$.
\end{theorem}

Theorem~\ref{interpolation-intro} states that $H_l(\beta)$ as a function of~$l$ from~$\Nb$ to~$\Qb(\beta)$ is  \emph{continuous} with respect to the $2$-adic topologies on~$\Nb$ and $\Qb(\beta)$,  and can even be extended to a~$2$-adic \emph{analytic} function in a $2$-adic neighborhood of~$\Nb$. The latter fact alone combined  with the recursion satisfied by Legendre polynomials already suffices to prove the finiteness of the set considered in Theorem~\ref{main-theorem-intro} (see Theorem~\ref{thm:finite-mult}). Further, to obtain the quantitative bound, we explicitly describe the radii of the disks on which $\check{H}_\beta$ is analytic in terms of the valuation of $\beta$. This interpolation result, combined with analysis of the local Newton polygons, yields \Cref{main-theorem-intro}. We note here that the analysis of Newton polygons for this function is not immediate since the coefficients of relevant power series belong to $\Q_2(\beta)$, which is hard to control. Luckily, we end up relying only on the total ramification of the extension $\Q_2(\beta)/\Q_2$: the coefficients of the relevant power series are naturally written as sums of terms, some of which have the same valuation. It is a special feature of totally ramified extensions of $\Q_2$ that a sum of an odd number of terms of valuation $v$ itself has valuation $v$, while a sum of an even number of terms of valuation~$v$ has valuation strictly bigger than~$v$;   allowing us to obtain an explicit bound in the end.

Finally, we show that this approach, while giving a finiteness result for every specific root, does not give a uniform bound for the number of times a given root can be repeated. We show
\begin{theorem}\label{counterexample-intro}
    For suitable choices of roots $\beta$ of Holt polynomials, the function $\check{H}_\beta(l)$, which is analytic in a certain neighborhood $N_h(\Z_2)$ of $\Z_2$, can have an arbitrary high number of roots on $N_h(\Z_2)$.
\end{theorem}

\subsection{Organization of the paper}
We introduce Legendre and Holt polynomials in Section~\ref{sec:holt}.
To render  the text accessible to a wider audience we have included in~\S\ref{sec:padic-analysis} a short summary of the results from $p$-adic analysis that we need. After proving the local analyticity result in~\S\ref{sec:analyticity}  we first present in~\S\ref{sec:finiteness} a ``soft'' proof of a simpler version of \Cref{main-theorem-intro}, which just asserts the finiteness of the set of $l$'s, and then, in~\S\ref{sec:NP} we give the bound on the size of this set, which requires a more technical $2$-adic calculation.
In~\S\ref{sec:counter-examples} we construct examples to show that an extension of Stieltjes's conjecture to the realm of~$2$-adic numbers is not possible (see Theorem~\ref{counterexample-intro}).
\vspace{2ex}

\paragraph{\textbf{Remark.}}
    As we were preparing the final version of this manuscript, the preprint \cite{deng2026irreducible} appeared, claiming to prove a strong form of Stieltjes conjecture that Legendre polynomials $P_l$ are irreducible (after removing the trivial zero at $x=0$) by a method different to the one used in the present paper. Several talks on our results were given by the second author in the period November 2025--February 2026 (at Hebrew University, CMS winter meeting, Newton Institute, IAS, and MIT).

\section*{Acknowledgments}
D.M. is very grateful to Akshay Venkatesh for his interest and encouragement, in particular for pointing out the connection of the work~\cite{mangoubi-weller-weiser2019harmonic} to Skolem's method. We thank Gal Binyamini, Avner Kiro, Jonathan Pila and Daniel Rosenblatt for useful discussions.
B.K. was supported by the Israel Science Foundation (ISF), Grant No.~2112/26.
D.M. was supported by the Israel Science Foundation (ISF), Grant No.~1219/25.
\section{Legendre polynomials and Holt's transformation}
\label{sec:holt}
\subsection{Definition and recursion}
The sequence of Legendre polynomials $(P_l)_{l=1}^{\infty}$ is a classical family of orthogonal polynomials appearing naturally in many problems of Mathematical Physics and classical analysis. They are defined as the unique family of polynomials of one variable such that $P_l$ is of degree~$l$ for which
\[ \int_{-1}^1 P_k(x)P_l(x)\, dx = 0 \text{ for }k\neq l,\]
and $P_k(1)=1$ (one can show that all zeros of $P_k$ lie in~$(-1, 1)$ \cite[\S3.3]{szego1975orthogonal}).
  It is not difficult to show Rodrigues' formula~\cite[\S4.3]{szego1975orthogonal}
  \[P_l(x):=\frac{1}{2^l l!}\frac{d^l}{dx^l}(x^2-1)^l\ .\]
\begin{lemma}[Legendre recursion]
\label{lem:legendre-recursion}
Set $\tilde{P}_l:=l! P_l$. 
    The following recursion is satisfied
    \[\tilde{P}_{l+1}(x)=(2l+1)x \tilde{P}_{l}(x)-l^2 \tilde{P}_{l-1}(x)\ .\]
\end{lemma}
\begin{proof}
 See, for example, \cite[\S4.5]{szego1975orthogonal}.

\end{proof}
If we roll out the recursion we get 
\begin{lemma}[Legendre recursion rolled out]
\label{lem:legendre-recursion-rolled}
For all $n\in\N$ and $l\in\N_0$
    \[\tilde{P}_{l+n}(x)=A_{n}(x, l)\tilde{P}_l(x)-l^2A_{n-1}(x, l+1)\tilde{P}_{l-1}(x)\]
    where $(A_{n}(x, l))_{n=0}^{\infty}$  is a sequence of  polynomials in $\Z[x, l]$, which satisfies the recursion
    \begin{equation}
    \label{eqn:A-rec}
    A_{n+1}(x, l)=(2l+1)xA_{n}(x, l+1)-(l+1)^2A_{n-1}(x, l+2),
    \end{equation}
    with initial conditions $A_1(x, l)=(2l+1)x$ and $A_0(x,l)=1$.
\end{lemma}
\begin{proof}
      This is a direct corollary of \Cref{lem:legendre-recursion}; proof is by induction on~$n$.
\end{proof}
\begin{remark}
     The highest degree monomial in $A_n$ is~$2^nx^nl^n$. The polynomial $A_n$ is of degree~$n$ with respect to~$x$ as well as with respect to~$l$.  
\end{remark}
Finally, we will use the following observation on the coefficients of $A_n(x, l) \in \Z[x,l]$.

\begin{lemma}\label{lem:chebyshev-in-legendre}
    Let $B_{n,k}(x) \in \Z[x]$ be polynomials such that $A_n(x,l)=\sum_{k=0}^n B_{n,k}(x)l^k$. Then $B_{n,n}(x)=U_n(x)$ is the Chebyshev polynomial of the second kind.
\end{lemma}
\begin{proof}
    By \Cref{lem:legendre-recursion-rolled}, specifically equation \eqref{eqn:A-rec}, comparing the $l^{n+1}$ terms on both sides we obtain the recursion
    $B_{n+1, n+1}(x)=2xB_{n,n}(x)-B_{n-1, n-1}(x)$. This recursion, combined with the initial conditions $B_{0,0}(x)=1$ and $B_{1,1}(x)=2x$, exactly describes the Chebyshev polynomials of the second kind (see \cite[eqs.~(4.1.7) and~(4.5.1)]{szego1975orthogonal}).
\end{proof}

\subsection{Integral representation and Holt's polynomials}
The first Laplace integral of Legendre polynomials is the following formula (see~\cite[\S4.8]{szego1975orthogonal}).
\begin{equation}
  \label{formula:laplace}
P_l(x):=\frac{1}{\pi}\int_0^{\pi} (x+\sqrt{x^2-1}\cos\alpha)^l \,d\alpha=
\sum_{k=0}^{\lfloor l/2 \rfloor} \binom{l}{2k}\binom{2k}{k} 4^{-k}x^{l-2k}(x^2-1)^{k}.
\end{equation}

  Holt~\cite{holt1913legendre}, while working on Stieltjes's conjecture, was the first to observe that the polynomials associated with this integral representation, have very nice $2$-adic properties, which can be used to study the conjecture.
  Define the polynomial~$H_{l, 0}$ of degree~$l$, corresponding to~$P_{2l}$ by
\begin{equation}
\label{def:Hl0}
     H_{l, 0}(z):=\frac{1}{\binom{2l}{l}}(z-4)^{l}
  P_{2l}\Bigl(\bigl(\frac{z}{z-4}\bigr)^{\frac{1}{2}}\Bigr).
\end{equation}
  
The polynomial $H_{l, 0}$ is normalized so that $H_{l, 0}(0)=1$. Observe that
\[P_{2l}(\cos\theta) =(-1)^l \binom{2l}{l}4^{-l}(\sin\theta)^{2l} H_{l, 0}(-4\cot^2\theta)\ .\]
 As the roots of $P_{2l}$ all lie in $[-1, 1]$, the roots of $H_{l, 0}$ lie in $(-\infty, 0)$.
Using the integral representation~\eqref{formula:laplace} an explicit expression for~$H_{l, 0}$ can be readily computed.
\begin{equation}
\label{formula:Hl0}
    H_{l, 0}(z)=\frac{1}{\binom{2l}{l}}\sum_{k=0}^l \binom{2l}{2k}\binom{2k}{k}z^{l-k}=\sum_{k=0}^{l} \frac{1}{\binom{2k}{k}}\binom{l}{k}^2 z^k\ .
\end{equation} 
Similarly, corresponding to $P_{2l+1}$,  define 
 
\begin{align}
\label{def:Hl1}
    H_{l, 1}(z):= \frac{1}{(2l+1)\binom{2l}{l}} (z-4)^l \bigl(\frac{z}{z-4}\bigr)^{-\frac{1}{2}} P_{2l+1}\Bigl(\bigl(\frac{z}{z-4}\bigr)^{\frac{1}{2}}\Bigr).
\end{align}
   Then, we have
   \[P_{2l+1}(\cos\theta)=(-1)^l(2l+1)\binom{2l}{l}4^{-l}(\sin\theta)^{2l}\cos\theta\cdot H_{l, 1}(-4\cot^2\theta)\]
   and
\begin{equation}
\label{formula:Hl1}
   H_{l, 1}(z)=\frac{1}{(2l+1)\binom{2l}{l}}\sum_{k=0}^l \binom{2l+1}{2k}\binom{2k}{k}z^{l-k}=\sum_{k=0}^l \frac{1}{(2k+1)\binom{2k}{k}}\binom{l}{k}^2 z^k\ .
\end{equation}

Throughout, we will use the notation $H_{l,\eps}$ to denote either one of $H_{l,0}$ or $H_{l,1}$ (so $\eps$ is only allowed to take values $0$ or $1$.)
\begin{remark}
    Note that the Legendre polynomials $P_{2l}(x), P_{2l+1}(x)$ are even and odd respectively. The change of variables $x^2=z/(z-4)$ gives a correspondence between pairs $\pm\alpha$ of \emph{nonzero} roots of Legendre polynomials with roots $\beta=4\alpha^2/(\alpha^2-1)$ of the corresponding Holt polynomial. Therefore Conjecture \ref{conj:stieltjes} is equivalent to the statement that no two distinct Holt polynomials can share a root. 
\end{remark}
\section{Basic \texorpdfstring{$p$}{p}-adic analysis}
\label{sec:padic-analysis}
\subsection{Analytic functions}
We briefly review some main principles of $p$-adic analysis. Our argument for proving \Cref{main-theorem-intro} is similar to the classical proof of the Skolem--Mahler--Lech theorem via $p$-adic interpolation (for an exposition see \cite{cassels1986local}). The main point is that for a fixed root $z$ of a Legendre polynomial, the function $P_{l}(z)$, after rescaling, becomes a ``$2$-adic analytic function'' of $l$ on $\Z_2$. However, since $\Z_2$ is compact, the function $P_l(z)$ can have at most finitely many zeroes or has to vanish identically on a disk (which is disproved using a separate argument). We now give a low-level summary of the theory for a reader familiar with the $2$-adic integers $\Z_2$, but not with deeper $p$-adic analysis. For a more detailed exposition of the subject see \cite{koblitz1984padic}, and for a nice exposition of Amice's theorem below see \cite{colmez2010fonctions}.

First we recall the following feature of analysis on $\Z_p$. We can define an analytic function on $\Z_p$ to be any function that can be represented by a converging power series $f(z)=\sum_{i\geqslant 0} a_i z^i$. The resulting space of functions, however, is not complete with respect to the supremum norm. In fact, due to the theory of Mahler series \cite{mahler1958interpolation}, any continuous function on $\Z_p$ can be approximated by analytic functions in the sup norm. This is not a desirable behavior for a theory analogous to complex analysis, but note that it is exactly in parallel with the space of real analytic functions on the closed interval in $\R$. Thus within this analogy $\Z_p$ should be thought of as the unit disk in $\R$, and to develop complex analysis we need to study functions defined on a larger set. This prompts the introduction of the field of complex $p$-adic numbers.

Motivated by the analogy with $\R$ one is led to consider an algebraic closure $\overline{\Q}_p$ of $\Q_p$. However the arithmetic of $\Q_p$ is richer than that of $\R$, and the extension $\overline{\Q}_p/\Q_p$ is of infinite degree, making $\overline{\Q}_p$ not complete with respect to the $p$-adic norm. Finally, let $\C_p$ denote the completion of an algebraic closure of $\Q_p$. This is a topological field equipped with the unique valuation $v:\C_p \to \Q\cup\{\infty\}$ which extends the $p$-adic valuation on $\Q_p$, and the associated norm $|z|_p=p^{-v_p(z)}$. 

For $z_0 \in \C_p$ and $r \in \R$ let $D(z_0,r)$ denote the closed disk~$\{z\in\C_p|\ |z-z_0|_{p}\leqslant r\}$. A function $f: D(z_0, r) \to \C_p$ is called \emph{analytic} on $D(z_0,r)$ if it can be represented by a convergent power series $f(z)=\sum_{k=0}^{\infty} a_k (z-z_0)^k$. Analytic functions on a disk form a Banach algebra, called the Tate algebra, with respect to the sup norm \[\|f\|_{D(z_0,r)}=\sup\{|f(z)|\colon \ z\in D(z_0,r)\}.\]

This norm is known as \emph{the Tate norm}. It can be read off from the Taylor series\footnote{this is not immediately obvious} by the following formula:
\[\|f\|_{D(z_0, r)}=\sup_k |a_k|_p r^k.\] 

The following example explicitly demonstrates the difference between supremum norms on $\C_p$ and $\Q_p$.

\begin{example}
    Consider the series $f(z)=\sum_{k=0}^{\infty} 2^k \binom{z}{k}$ on the unit disk $D(0,1)\subset\C_2$. Each binomial $g_k(z)=\binom{z}{k}$ is an integer-valued function on $\Z_2$, and therefore $f(z)$ defines a continuous function on $\Z_2$. However, this function does not extend to an analytic function on the unit disk $D(0,1)$ in $\C_2$. The reason is that while $g_k$ takes small values on $\Z_2$, it takes large values at other points of $D(0,1)$. By the formula for the Tate norm, $\|g_k\|_{D(0, 1)}=1/|k!|_2$. The $2$-adic valuation of the factorial is $v_2(k!)=k-s_2(k)$, where $s_2$ is the sum of the base $2$ digits of $k$. This shows that the series defining $f(z)$ does not converge, since $\|2^k g_k\|=2^{-s_2(k)}$ does not tend to~$0$ as $k$ tends to~$\infty$.
\end{example}

The previous example is important for us, since Legendre polynomials can be represented as series with binomial terms. To sidestep the issue of binomials having large norms, we consider the algebra of $h$-\emph{locally} analytic functions $LA_h(\Z_p, K)$, following Amice \cite{amice1964interpolation}; see also \cite{colmez2010fonctions}, especially \cite[Section I.4]{colmez2010fonctions}. Note that in $p$-adic analysis ``locally analytic'', in the sense of being represented by a power series near every point, is not a well-behaved concept (in particular, locally-analytic functions on a disk are not necessarily analytic). We use the following notion instead. Let $K$ denote a finite extension of $\Q_p$. The space $LA_h(\Z_p, K)$ consists of functions $f: \Z_p \to K$ such that for every $z_0 \in \Z_p$ the function $f$ restricted to $\Z_p \cap D(z_0, p^{-h})$ can be written as a power series centered at $z_0$ with coefficients in $K$ which converges in $\C_p$ on $D(z_0, p^{-h})$. Such a power series representation is necessarily unique, and $f$ extends to a function on the $p^{-h}$-neighborhood of $\Z_p$:  \[N_{h}(\Z_p)\colonequals \left\{z \in \C_p: |z-t| \leqslant p^{-h} \text{ for some } t \in \Z_p\right\}.\]

Note that $h$ is a valuative parameter, so that the higher the $h$ the smaller the domain $N_h(\Z_p)$.

\begin{proposition}[Main properties of locally-analytic functions]\label{prop:analytic-functions-basics}
\noindent
\begin{enumerate}
    \item\label{prop-BanachA} The space $LA_h(\Z_p, K)$ is a Banach algebra over $K$ with respect to the sup norm $\| \cdot \|_{LA_h}$ on $N_h(\Z_p)$.
    \item\label{prop-Binom-norm} Let \[g_n(z):=\binom{z}{n}.\] We have \[\left\|g_n\right\|_{LA_h}=\left| \lfloor np^{-h}\rfloor! \right|_p^{-1}.\]
    \item\label{prop:finiteness} For any $f \in LA_h(\Z_p, K)$ the zero-set of $f$ on $\Z_p$ is a finite union of points and closed disks of radius $p^{-h}$. In particular, such a function has only finitely many zeroes on $\Z$ if it does not vanish identically on any residue class modulo $p^h$.
\end{enumerate}
\end{proposition}
\begin{proof}
Parts \eqref{prop-BanachA} and \eqref{prop-Binom-norm} follow from the work of Amice \cite{amice1964interpolation}; see also \cite[Section I.4]{colmez2010fonctions}. Part \eqref{prop:finiteness} follows from Strassman's theorem \cite[Theorem 4.1]{cassels1986local}: a nonzero analytic function on a disk has only finitely many zeroes in~$\Z_p$ \footnote{ Nonzero analytic functions whose Taylor coefficients belong to $K$ also have at most finitely many zero on the disk in $\C_p$, but this is a harder theorem, which can be proved using the theory of Newton polygons; see, for example, \cite[ Ch. IV\S4]{koblitz1984padic} }. A function $f \in LA_h(\Z_p, K)$ is analytic on the disks $D(a, p^{-h})$ for $a=0,1,\dots, p^h-1$, and $D(a, p^{-h})\cap \Z =\{n\in \Z\colon\ n\equiv a \mod p^h  \}.$
\end{proof}
\begin{remark}
    In fact, Amice \cite{amice1964interpolation} proves that the binomials $e_n(z)=\lfloor np^{-h}\rfloor! \binom{z}{n}$ form an orthonormal basis of $LA_h(\Z_p, K)$ in the sense that every element $f \in LA_h(\Z_p, K)$ can be written uniquely as a series $f=\sum_k a_k e_k$ with $a_k \in K$ tending to zero, and $\|f\| = \sup_k |a_k|$. We will not use this in what follows.
\end{remark}

\begin{remark}
    We have used the language of norms, as it makes clear the analogy with complex analysis. However, in practical computation, it is more convenient to work with valuations, which we will do for now on. The connection between the two, whenever it is a norm on $\Z_p, \C_p$ or some algebra of functions, is $|\bullet|=p^{-v(\bullet)}$. If the prime $p$ is not clear from the context, we indicate it in the subscript $|\bullet|_p, v_p(\bullet)$.
\end{remark}
\subsection{The Legendre--Kummer lemma}
The following standard elementary $p$-adic lemma on valuations of factorials and binomials will be used many times in calculations. 

\begin{lemma}[Legendre+Kummer]
    The $p$-adic valuation of a factorial is \[v_p(n!)=\frac{n-s_p(n)}{p-1} \leqslant n/(p-1),\] where $s_p$ is the base $p$ sum-of-digits function. In particular, for the $2$-adic valuation of a binomial we have 
    \[v_2\Biggl(\!\binom{n}{k}\!\Biggr)=s_2(k)+s_2(n-k)-s_2(n)\ .\]
    Specializing to the case $k=1$ we get
    \[v_2(n)=1+s_2(n-1)-s_2(n)\ ,\]
    Further,  for all $a, b\in\N$ \[s_2(a)+s_2(b)\geqslant s_2(a+b)\ .\]
\end{lemma}
\begin{proof}
    The claims about valuation can be found in many elementary number theory texts like \cite[Chapter 6]{roettger2025enchant} (the historical reference is \cite{legendre1830theorie}). The subadditivity of the sum-of-digits function $s_2$ follows from the inequality $0\leqslant v_2(\binom{a+b}{a})=s(a)+s(b)-s(a+b)$. 
\end{proof}

\subsection{Wahab's theorem}\label{sec:Wahab}
We will need information about the $2$-adic properties of roots of Legendre polynomials. The necessary result has already been obtained by Wahab \cite{wahab1952new}. The key to Wahab's result is the nice formula (easily obtained from the Rodrigues formula, see~\cite{grosswald1950}) for the shifted Legendre polynomial 
\[P_l(1+2u)=\sum_{k=0}^l\binom{l+k}{l}\binom{l}{k}u^k,\]
which allows for a direct analysis of the $2$-adic Newton polygon of the shift, and thus gives the valuations of shifted roots of $P_l$. 

\begin{theorem}[{\cite[Th.~3.1]{wahab1952new}}]\label{theorem:Wahab}
    Suppose $\alpha\in\Qbar\subset \Qbar_2$ is a nonzero root of some Legendre polynomial $P_l$. Suppose $l=2^{e_1}+2^{e_2}+\dots +2^{e_s}$ is the base $2$ expansion of $l$. Then for some $i$ we have
    \[v_2(\alpha-1)=1-2^{-e_i}.\]
    Moreover, for each $i$, $P_l$ has exactly $2^{e_i}$ roots  in $\Qbar_2$ with $v_2(\alpha-1)=1-2^{-e_i}$; in particular, every nonzero root $\alpha$ satisfies $v_2(\alpha-1)>0$.
\end{theorem}
Changing the variables, this gives the following parallel result for the Holt polynomial.

\begin{corollary}\label{cor:Wahab-for-H}
 Suppose $H_{l, \eps}(\beta)=0$ and $\beta\neq 0$. Suppose $l=2^{e_1}+\dots +2^{e_s}$ is the base $2$ expansion of $l$. Then $v_2(\beta)=2^{-e_i}$.
\end{corollary}
\begin{proof} 
Recall that $\beta$ can be written in terms of the corresponding root $\alpha$ of a the Legendre polynomial $P_{2l+\eps}$ as follows:
\[\beta=\frac{4\alpha^2}{\alpha^2-1}= \frac{4(\alpha-1)^2+8(\alpha-1)+4}{(\alpha-1)(1+\alpha)}.\]
Using Wahab's \Cref{theorem:Wahab} we can compute the valuation of the right hand side. Since $v_2(\alpha-1)>0$ (for $\alpha \neq 0$) we have
\[v_2(4\alpha^2)=v_2(4(\alpha-1)^2+8(\alpha-1)+4) =2, \]
and since $v_2(\alpha-1)<1$ we have
\[v_2(1+\alpha)=v_2(2+(\alpha-1))=v_2(\alpha-1)\ .\]
We conclude that
\begin{equation}
\label{eqn:valphavbeta}
v_2(\beta)=2-2v_2(\alpha-1).
\end{equation}
The result now follows from \Cref{theorem:Wahab} applied to $P_{2l+\epsilon}$.
\end{proof}

\subsection{Auxiliary results}
 Finally, we will use two simple statements whose proofs require a bit more theory than what we have introduced. 

 \begin{proposition}\label{prop:odd-even}
     Let $\alpha$ be a root of a Legendre polynomial, and consider the field $K=\Q_2(\alpha)$. Suppose $x_1, \dots, x_n\in K$ are elements of the same valuation. Then
     \[\begin{cases}
     v(x_1+\dots+x_n)=v(x_1), \text{ if } n \text{ is odd,}\\
     v(x_1+\dots + x_n)>v(x_1), \text{ if } n \text{ is even.}
     \end{cases}\]
 \end{proposition}
 \begin{proof}
     By Wahab's theorem, the field $K$ is a totally ramified extension of $\Q_2$. The result follows, since the residue field of $K$ is $\F_2$. 
 \end{proof}

 \begin{proposition}\label{prop:single-root}
     Suppose $f(z)=a_0+a_1z+\dots$ is a power series with coefficients in a finite extension $K$ of $\Q_p$ which converges on the closed disk $D(0,1)$. Suppose that for $i>1$, $v(a_i)>v(a_1)$. Then $f(z)$ has at most one root in $D(0,1)$, and that root belongs to $K$. If additionally $v(a_0)\geqslant v(a_1)$, then $f(z)$ has exactly one root.
 \end{proposition}
 \begin{proof}
     This follows from the theory of Newton polygons for power series (which, more generally, allows one to describe the roots of a power series and their valuations); see, for example, \cite[Chapter IV.4]{koblitz1984padic}. The result we need is an immediate corollary of \cite[Theorem IV.4.14]{koblitz1984padic}.   
 \end{proof}
\section{\texorpdfstring{$2$}{2}-adic interpolation of Holt's polynomials}
\label{sec:analyticity}

We prove the existence of the $2$-adic interpolation of Holt polynomials $H_l(\beta)$ as functions of $l$, whenever $v_2(\beta)>0$.

\begin{theorem}\label{thm:holt-interpolations}
    Suppose $\beta \in \overline{\Q}_2$ is a constant for which $v_2(\beta)>0$. Consider $\check{H}_{\beta, \eps}(l):=H_{l, \eps}(\beta)$ as a function from $\N$ to $\Q_2(\beta)$. There exists $h\in\N$ such that $\check{H}_{\beta, \eps}$ extends to a function $\check{H}_{\beta, \eps}\in LA_h(\Zb_2, \Qb_2(\beta))$. Moreover, one can take $h=\lceil 2-\log_2(v_2(\beta))\rceil$. In particular, if $\beta$ is a root of $H_{l_0,\eps}$, then one can take $h=2+\lceil \log_2 l_0 \rceil$.
\end{theorem}
\begin{proof}
    By Amice's theorem the polynomial function 
    \[g_k(l)=\binom{l}{k}\] satisfies
    \[v_{LA_h}(g_k)=-v_2(\lfloor k/2^h\rfloor!).\]
  Consider  the following two  series
    \[\check{H}_{\beta, 0}=\sum_{k=0}^{\infty} \frac{\beta^k}{\binom{2k}{k}}g_k^2, \text{ and } \check{H}_{\beta, 1}=\sum_{k=0}^\infty \frac{\beta^k}{(2k+1)\binom{2k}{k}}{g_k^2}.\]
    For positive integers $l$ we have $\check{H}_{\beta, \eps}(l)=\tilde{H}_{\beta, \eps}(l)$ by formulas \eqref{formula:Hl0} and \eqref{formula:Hl1}. We claim that these series converge in  $LA_h(\Zb_2, \Qb_2(\beta))$ if $h$ is large enough. Indeed, for the first series we compute
\begin{align}
\label{eqn:linear-growth}
v_{LA_h}\left(\frac{\beta^k}{\binom{2k}{k}}{g_k^2}\right)=kv_2(\beta)- s_2(k)-2v_2(\lfloor k/2^h\rfloor!)=k v_2(\beta) -s_2(k)-2\lfloor k/2^h \rfloor +2s_2(\lfloor k/2^h\rfloor) 
\end{align}
    The latter expression tends to~$\infty$ as $k\to\infty$ as long as $v_2(\beta)>2^{1-h}$. We can choose~$h$ satisfying this inequality since $v_2(\beta)>0$. The same conclusion holds for the other series, as each term in the series for $\check{H}_{\beta,1}$ differs from the corresponding term in the series of $\check{H}_{\beta,1}$ by an odd factor.

    Finally, the last assertion follows: by Corollary \ref{cor:Wahab-for-H} of Wahab's theorem, if $\beta$ is a root of $H_{l,\eps}$, then $v_2(\beta)>0$.
\end{proof}

\section{Finite multiplicity}
\label{sec:finiteness}
In this section we prove the following weaker version of \Cref{main-theorem-intro}.

\begin{theorem}\label{thm:finite-mult}
    Suppose $\alpha\neq 0$ is a root of the Legendre polynomial $P_l$. Then the set 
    \[\{n\in \N \colon P_n(\alpha)=0\}\] is finite.
\end{theorem}

Let $z_0\in\Zb_2$ be a point. By part \eqref{prop:finiteness} of \Cref{prop:analytic-functions-basics}, a nonzero function in~$LA_{h}(\Zb_2, K)$ has either a finite number of zeros in the disk  $D(z_0, p^{-h})$, or it vanishes  in that whole  disk. Since~$\Zb_2$ is covered by a finite number of such disks, it only remains to verify that $\check{H}_{\beta, 0}$ and $\check{H}_{\beta, 1}$ cannot  vanish identically on such disks. We will do so by using the recursion formula for Legendre polynomials. We first prove that these recursions are nondegenerate.

\begin{lemma}\label{lem:A-not-zero}
     For any $n\in \N$ and a nonzero Legendre root~$\alpha\in \overline{\Qb}$ the polynomial $A_n(\alpha, l)\in \overline{\Qb}[l]$ is not the zero polynomial.
\end{lemma}

\begin{proof}
Recall from Lemma \ref{lem:chebyshev-in-legendre} that we can write $A_n(x, l)=\sum_{k=0}^n B_{n,k}(x) l^k$, where $B_{n,n}(x)=U_n(x)$ is the Chebyshev polynomial of the second kind. Thus it suffices to show that if $\alpha$ is a root of a Legendre polynomial, then $U_n(\alpha)\neq 0$.

An explicit formula for Chebyshev polynomials of the second kind is
    \[ U_n(\alpha)= \sum_{k=0}^{\lfloor n/2 \rfloor} \alpha^{n-2k} (\alpha^2-1)^{k} \binom{n+1}{2k+1}\ . \]
    We claim that for $1\leqslant k \leqslant \lfloor n/2 \rfloor$ we have 
    \begin{equation}
    \label{ineq:chebyschev-2-powers}
    v_2\Biggl(\binom{n+1}{2k+1}\Biggr)\geqslant v_2(n+1).
    \end{equation}
    Indeed, 
    
    \begin{multline*}
    v_2\Biggl(\binom{n+1}{2k+1}\Biggr) = s_2(2k+1)+s_2(n-2k)-s_2(n+1)=s_2(2k)+1+s_2(n-2k)-s_2(n+1) \\ \geqslant s_2(n)+1-s_2(n+1)=v_2(n+1).
    \end{multline*}
    If $\alpha\neq0$ is a root of $P_l$, then by Wahab's theorem \ref{theorem:Wahab} we have $v_2(\alpha)=0$ and $v_2(\alpha-1)>0$, and so $v_2(U_n(\alpha))=v_2(n+1)$. In particular, $U_n(\alpha)\neq 0$.\qedhere

\end{proof}
We are now ready to prove the main result of this section.

\begin{proof}[Proof of Theorem~\ref{thm:finite-mult}]
It is enough to prove the same result for the Holt polynomials $H_{l,\eps}$. Suppose $\beta$ is a root of a Holt polynomial (so that $\alpha=\left(\beta/(\beta-4)\right)^{1/2}$ is a root of a Legendre polynomial). Then \Cref{thm:holt-interpolations} shows that there are locally analytic functions $\check{H}_{\beta,\eps}(l)$ such that for every positive integer $l$ we have $\check{H}_{\beta,\eps}(l)=H_{l, \eps}(\beta)$. To conclude we just need to show that the functions $\check{H}_{\beta, \eps}$ cannot identically vanish on a residue disk.

Assume that $\check{H}_{\beta, \eps}=0$ on $D(x, 2^{1-h})$ for some $x\in \Zb_2$ and some $h \in \N$. There exists $m\in\N_0$, $m<2^{h-1}$, such that $\tilde{H}_{\beta, \eps}(m+k 2^{h-1})=0$ for all $k\in\N_0$, or equivalently, $P_{2m+\eps+ k 2^{h}}(\alpha)=0$ for some $\alpha\in\Qbar$ and all $k\in\N_0$. Let $\tilde{P}_l(x)=l! P_l(x)$. We know that for all $l\in\N$ such that $l\equiv 2m+\eps \mod{2^h}$ we have $\tilde{P}_l(\alpha)=0$, and by Lemma~\ref{lem:legendre-recursion-rolled} we have
\[\tilde{P}_{l+2^h}(\alpha)=A_{2^h}(\alpha, l)\tilde{P}_{l}(\alpha)-l^2  A_{2^h-1}(\alpha, l+1)\tilde{P}_{l-1}(\alpha).\]
By Lemma~\ref{lem:A-not-zero}, we can find $l_0\in\N$,  such that $l_0\equiv 2m+\eps\mod{2^h}$ and
$A_{2^h-1}(\alpha, l_0+1)\neq 0$.
It then follows from the above recurrence that $\tilde{P}_{l_0}(\alpha)=\tilde{P}_{l_0-1}(\alpha)=0$.
Finally, it follows from the standard Legendre recursion, Lemma~\ref{lem:legendre-recursion}, that $\tilde{P}_0(\alpha)=0$, which is clearly false.
\end{proof}
\section{Bounding the number of roots}
\label{sec:NP}
In  section~\ref{sec:analyticity} we have seen that if $\beta\in\Cb_2$ is such that $v(\beta)=2^{-n}$ then $\check{H}_{\beta, \epsilon}(l):=H_{l, \eps}(\beta)$
is analytic on any disk $B(l_0, 2^{-n-2})$ where $l_0\in\Zb$.
Note that $\Zb$ can be covered by $2^{n+2}$ such disks.
Due to Wahab's theorem, in any such disk, if $\lfloor l_0/2^n\rfloor$ is even, then $\check{H}_{\beta, \epsilon}$ has no roots lying in~$\Zb$. So we are left with $2^{n+1}$ relevant disks. In this section we will prove the following result.
\begin{theorem}\label{thm:unique-in-disk}
Let $\beta\in\Qbar_2$ be a root of of a Holt polynomial with valuation $v(\beta)=2^{-n}$. Assume $l_0\in\Zb$ is such that $\lfloor l_0 2^{-n}\rfloor$ is odd. Then, the function $\check{H}_{\beta, \epsilon}$ has at most one root in the disk $B(l_0, 2^{-n-2})$.
\end{theorem}
This immediately implies our main theorem.
\begin{corollary}[\Cref{main-theorem-intro}]
Suppose $\alpha\neq 0$ is a root of the Legendre polynomial $P_l$. Then the set 
\[\{n\in \N \colon P_n(\alpha)=0\}\]
has cardinality at most $2l$.
\end{corollary}
\begin{proof}
The root $\alpha$ of $P_n$ corresponds to a root $\beta$ of the Holt polynomial $H_{\lfloor l/2 \rfloor, \epsilon}$, $\eps$ depending on the parity of~$n$. The valuation of $\beta$ is thus $2^{-m}$ with $2^m\leqslant l/2$ by Wahab's theorem, \Cref{cor:Wahab-for-H}. On each of the $2^{m+1}$ relevant residue disks modulo $2^{m+2}$ the functions $H_{\beta, 0}$ and $H_{\beta, 1}$ have at most one root. Altogether this gives at most $4\cdot 2^{m} \leqslant 2l$ elements in $\{n\in \N \colon P_n(\alpha)=0\}.$    
\end{proof}

To prove \Cref{thm:unique-in-disk} we will appeal to \Cref{prop:single-root}, which itself is a specific instance of the theory of Newton polygons. We thus will need to estimate the valuations of coefficients of Taylor expansions of $\check{H}_{\beta, \epsilon}$. It turns out that there is no difference in the calculation for $\check{H}_{\beta,0}$ and $\check{H}_{\beta, 1}$, since the formulas for the relevant series, \eqref{formula:Hl0}, \eqref{formula:Hl1}, differ by odd integer factor added to the terms of the series. To lower the notational burden, we will write the calculations for the case of $\check{H}_\beta\colonequals \check{H}_{\beta,0}$. The other case can be obtained verbatim by simply inserting the missing odd factors in the terms of the series.

Finally, the coefficients of relevant power series will be sums of terms whose valuations are expressions involving sums of digits and valuation functions. The combinatorial problem of minimizing such expressions will be postponed to \Cref{sec:auxiliary} to clarify the logic of the argument.

\subsection{Computing Newton Polygons}
For the rest of this section we are using the setting and notation of \Cref{thm:unique-in-disk}.

It will be convenient to scale and shift $\check{H}_{\beta}$ to the unit disk. Let 
\begin{align}
\label{eqn:scale&shift}
    F(x):=\check{H}_{\beta}(l_0+2^{n+2}x)=\sum_{k=0}^{\infty} \frac{\beta^k}{\binom{2k}{k}}\binom{l_0+2^{n+2}x}{k}^2.
\end{align}
The function $F$ is analytic in $\Ol_{\Cb_2}$.
We consider its Taylor expansion 
\[F(x)\equalscolon \sum_{k=0}^{\infty} a_k x^k,\]

\subsubsection{Valuation of~$a_1$}
\begin{proposition}[Computing $v(a_1)$]
\label{prop:va1}
        \[v(a_1)=\begin{cases}
    3+l_0 2^{-n}-s(l_0), & \ 2^n\leqslant l_0<2^{n+1}.\\
    2+l_02^{-n}-s(l_0), &  2^{n+1}+2^n\leqslant l_0 <2^{n+2}.
    \end{cases}\]
\end{proposition}

\begin{proof}
Directly expanding each term in \eqref{eqn:scale&shift}, and collecting linear terms gives
\[a_1=2^{n+3}\sum_{k=1}^{l_0}\frac{\beta^k}{\binom{2k}{k}}\binom{l_0}{k}^2 \sum_{j=l_0-k+1}^{l_0}\frac{1}{j}=:\sum_{k=1}^{l_0} a_{1, k}.\]

We observe that on the interval $l_0-k+1\leqslant j\leqslant l_0$ there is exactly one~$j$ of maximal valuation.
Hence,
\[v(a_{1, k})=n+3+k2^{-n}+v\left(\!\binom{l_0}{k}\!\right)+s(l_0-k)-s(l_0)-\max_{l_0-k+1\leqslant j\leqslant l_0} v(j)\ .\]
In particular,
\[v(a_{1, l_0})=n+3+l_02^{-n}-s(l_0)-\max_{1\leqslant j\leqslant l_0} v(j)=\begin{cases}
3+l_02^{-n}-s(l_0),& \quad\quad \quad 2^n\leqslant l_0<2^{n+1},\\
2+l_02^{-n}-s(l_0),& 2^{n+1}+2^n\leqslant l_0<2^{n+2}.
\end{cases}\]
We will show that, in fact, $v(a_1)=v(a_{1, l_0})$.

Assume first that $2^n\leqslant l_0<2^{n+1}$. Then 
\begin{align}
\label{eqn:a1k-a1l0-caseA}
v(a_{1, k})- v(a_{1, l_0})=-(l_0-k)2^{-n}+v\left(\!\binom{l_0}{k}\!\right)+s(l_0-k)+n-\max_{l_0-k+1\leqslant j\leqslant l_0} v(j).
\end{align}

\begin{itemize}
\item If $s(l_0-k)\geqslant 2$, then the RHS of~\ref{eqn:a1k-a1l0-caseA} is positive,
since $v(j)\leqslant n$ for any~$j\leqslant l_0$, and $(l_0-k)2^{-n}< 2$.

\item If $s(l_0-k)=1$, then $l_0-k=2^m$ for some $0\leqslant m\leqslant n$. If $m<n$ the RHS of~\eqref{eqn:a1k-a1l0-caseA} is positive since $v(j)\leqslant n$ for any~$j\leqslant l_0$, and $(l_0-k)2^{-n}<1$.
If $m=n$ then $v(j)\leqslant n-1$ for any $j\geqslant l_0-k+1=2^{n}+1$ and we get again that the RHS of~\eqref{eqn:a1k-a1l0-caseA} is positive.
\end{itemize}
We have shown that  $v(a_{1, k})\geqslant v(a_{1, l_0})$ with equality if and only if $k=l_0$ and we conclude that 
$v(a_1)=v(a_{1, l_0})$ as desired.

Now assume that $2^{n+1}+2^n\leqslant l_0<2^{n+2}$. In this case we compute
\begin{align}
\label{eqn:a1k-a1l0-caseB}
    v(a_{1, k})-v(a_{1, l_0})& = -(l_0-k) 2^{-n}+v\left(\!\binom{l_0}{k}\!\right)+s(l_0-k)+n+1-\max_{l_0-k+1\leqslant j\leqslant l_0} v(j).
\end{align}
\begin{itemize}
    \item In case $s(l_0-k)\geqslant 2$
    \begin{itemize}
        \item if $1\leqslant l_0-k<2^{n+1}$ then $(l_0-k)2^{-n}<2$ and  clearly, the RHS of~\eqref{eqn:a1k-a1l0-caseB} is positive.
        \item if $2^{n+1}\leqslant l_0-k<2^{n+1}+2^n$ then $(l_0-k)2^{-n}<3$ and $\max_{l_0-k+1\leqslant j\leqslant l_0} v(j)=n$. So  the RHS of~\eqref{eqn:a1k-a1l0-caseB} is positive.
        \item if $2^{n+1}+2^n\leqslant l_0-k<2^{n+2}$ then $(l_0-k)2^{-n}<4$ and $\max_{l_0-k+1\leqslant j\leqslant l_0} v(j)\leqslant n-1$. So  the RHS of~\eqref{eqn:a1k-a1l0-caseB} is positive.
    \end{itemize}
    \item In case $s(l_0-k)=1$ then $(l_0-k)=2^{m}$ for some $0\leqslant m\leqslant n+1$. If $m<n$
    then the RHS of~\eqref{eqn:a1k-a1l0-caseB} is positive. If $m=n$ or $m=n+1$ the
    RHS of~\eqref{eqn:a1k-a1l0-caseB} is~$0$.
\end{itemize}

We have shown that $v(a_{1, k})\geqslant v(a_{1, l_0})$ with equality if and only if $k=l_0$ or $k=l_0-2^n$ or $k=l_0-2^{n+1}$. By \Cref{prop:odd-even}, we conclude that 
$v(a_1)=v(a_{1, l_0})$ as desired.
\end{proof}
\subsubsection{Valuation of~$a_2$}
In the next step of our investigation of the Newton Polygon of~$F$ we show that $v(a_2)$ is strictly bigger than~$v(a_1)$.
\begin{proposition}
    \[v(a_2)>v(a_1).\]
\end{proposition}
\begin{proof}

From expanding each term in~\eqref{eqn:scale&shift} it follows that we can express~$a_2$ as the following infinite convergent series: $a_2=\sum_{k=1}^{\infty} a_{2, k}$ 
where for $k\leqslant l_0$
\[a_{2, k} = 2^{2n+4}\frac{\beta^{k}}{\binom{2k}{k}}\binom{l_0}{k}^2\sum_{J\in\binom{S_k}{2}}\prod_{j\in J}\frac{1}{j}\]
and where
$S_k$ is the multiset
\[S_k=[l_0, l_0, l_0-1, l_0-1, \ldots, l_0-k+1, l_0-k+1],\]
and  for $k>l_0$ we have
\[a_{2, k}=2^{2n+4}\frac{\beta^k}{\binom{2k}{k}}\frac{1}{\binom{k}{l_0}^2(k-l_0)^2}.\]

Assume first that $2^n\leqslant l_0<2^{n+1}$. 
\begin{itemize}
    \item If $1\leqslant k\leqslant l_0$ we have
\begin{align*}
v(a_{2, k})-v(a_1)= -(l_0-k)2^{-n}+v\left(\!\binom{l_0}{k}\!\right)+s(l_0-k)+2n+1-2\max_{l_0-k+1\leqslant j\leqslant l_0} v(j).
\end{align*}
Since $v(j)\leqslant n$ for all~$j$ and $(l_0-k)2^{-n}<2$ we may conclude that
\begin{align}
\label{eqn:a2k-a1-caseA1}
    v(a_{2, k})-v(a_1)>0.
\end{align}
\item
If $k>l_0$ we have
\begin{align*}
    v(a_{2, k})-v(a_1)&= (k-l_0)2^{-n}+s(k)-2v(k-l_0)-2s(k-l_0)-s(l_0)+2n+1
\end{align*}
From the auxiliary Lemma~\ref{lem:G1-properties} it follows that
\begin{align}
\label{eqn:a2k-a1-caseA2}
v(a_{2, k})-v(a_1)\geqslant 0
\end{align}
 and 
\begin{equation}
\label{eqn:a2k-a1=0-caseA2}
v(a_{2, k})-v(a_1)= 0 \text{ iff } k-l_0=2^n \text{ or } k-l_0=2\cdot2^n \text{ or } k-l_0=3\cdot 2^n \text{ or } k-l_0=4\cdot2^n.
\end{equation}

\end{itemize}
From~\eqref{eqn:a2k-a1-caseA1},~\eqref{eqn:a2k-a1-caseA2} and~\eqref{eqn:a2k-a1=0-caseA2} we conclude (by \Cref{prop:odd-even}) that $v(a_2)-v(a_1)>0$ as desired.

Assume now that $2^{n+1}+2^n\leqslant l_0<2^{n+2}$. If $1\leqslant k\leqslant l_0$ we have
\begin{align*}
v(a_{2, k})-v(a_1)&=-(l_0-k)2^{-n}+v\left(\!\binom{l_0}{k}\!\right)+s(l_0-k)+2n+2-2\max_{l_0-k+1\leqslant j\leqslant l_0} v(j) 
\end{align*}

\begin{itemize}
\item If $1\leqslant l_0-k< 2^{n+1}$
then
\begin{equation}
\label{eqn:a2k-a1-caseB1}
v(a_{2, k})-v(a_1)=-(l_0-k)2^{-n}+v\left(\!\binom{l_0}{k}\!\right)+s(l_0-k)
\end{equation}
If $s(l_0-k)\geqslant 2$ then the last expression is positive, since $(l_0-k)2^{-n}<2$. If $s(l_0-k)=1$ then $l_0-k=2^m$ for some $0\leqslant m\leqslant n$. If $m<n$ then  $(l_0-k)2^{-n}< 1$ and the RHS of~\eqref{eqn:a2k-a1-caseB1} is positive. If $m=n$ then we get that the RHS of~\eqref{eqn:a2k-a1-caseB1} is~$0$.
\item If $2^{n+1}+1\leqslant l_0-k<2^{n+2}$ then
\begin{align}
\label{eqn:a2k-a1-caseB2}
v(a_{2, k})-v(a_1)&\geqslant 2-(l_0-k)2^{-n}+v\left(\!\binom{l_0}{k}\!\right)+s(l_0-k).
\end{align}
 Necessarily,  $s(l_0-k)\geqslant 2$, and since $(l_0-k)2^{-n}<4$ it follows that  the RHS of~\ref{eqn:a2k-a1-caseB2} is positive.
 \item If $k>l_0$ then 
\begin{align*}
    v(a_{2, k})-v(a_1)&= (k-l_0)2^{-n}+s(k)-2v(k-l_0)-2s(k-l_0)-s(l_0)+2n+2.
\end{align*}
From Lemma~\ref{lem:G1-properties} it follows that
\begin{align*}
v(a_{2, k})-v(a_1)\geqslant 0
\end{align*}
 and 
 \begin{equation*}
v(a_{2, k})-v(a_1)= 0 \text{ iff } k=l_0+2^n \text{ or } k=l_0+2\cdot2^n. 
\end{equation*}
\end{itemize}
We have shown that $v(a_{2, k})-v(a_1)\geqslant 0$ for all $k\geqslant 1$ with equality if and only if $k=l_0-2^n$ or $k=l_0$ or $k=l_0+2^n$ or $k=l_0+2^{n+1}$. It follows (from \Cref{prop:odd-even}) that
$v(a_2)-v(a_1)>0$ as desired.
\end{proof}
\subsubsection{Valuation of~$a_p$ for $p\geqslant 3$}
\begin{proposition}
    For all $p\geqslant 3$ \[v(a_p)\geqslant v(a_1)+1.\]
\end{proposition}
\begin{proof}
From~\eqref{eqn:scale&shift} it follows that we can express~$a_p$ as an infinite convergent series: $a_p=\sum_{k=1}^{\infty} a_{p, k}$, where for $k\leqslant l_0$ we have
\[a_{p, k} = 2^{p(n+2)}\frac{\beta^{k}}{\binom{2k}{k}}\binom{l_0}{k}^2\sum_{J\in\binom{S_k}{p}}\prod_{j\in J}\frac{1}{j}\]
for the multiset $S_k$ 
\[S_k=[l_0, l_0, l_0-1, l_0-1, \ldots, l_0-k+1, l_0-k+1],\]
and for $k>l_0$ we have
\[a_{p, k}=2^{p(n+2)}\frac{\beta^k}{\binom{2k}{k}}\frac{1}{\binom{k}{l_0}^2(k-l_0)^2}\sum_{J\in \binom{T_k}{p-2}}\prod_{j\in J}\frac{1}{j}\]
where $T_k$ is the multiset
\[T_k=[l_0, l_0, l_0-1, l_0-1, \ldots, 1, 1, -1, -1, -2, -2, \ldots -(k-l_0-1), -(k-l_0-1)]. \]

Assume first that $2^n\leqslant l_0<2^{n+1}$. 
\begin{itemize}
    \item If $k\leqslant l_0$ we have
\begin{align*}
v(a_{p, k})&\geqslant p(n+2)+k2^{-n}+v\left(\!\binom{l_0}{k}\!\right)-s(l_0)+s(l_0-k)-pn\\ 
&=2p+k2^{-n}+v\left(\!\binom{l_0}{k}\!\right)-s(l_0)+s(l_0-k). 
\end{align*}
Hence,
\begin{align*}
  v(a_{p, k})-v(a_1)\geqslant 2p-3-(l_0-k)2^{-n}+v\left(\!\binom{l_0}{k}\!\right)+s(l_0-k)> 1.
\end{align*}
\item
If $k>l_0$
\begin{align*}
    v(a_{p, k})-v(a_1)&\geqslant (k-l_0)2^{-n}+s(k)-2v(k-l_0)-2s(k-l_0)-s(l_0)+2n+1\\
    &+(p-2)(n+2)-\max_{J\in \binom{T_k}{p-2}} \sum_{j\in J}v(j)
\end{align*}
The estimate $v(a_{p, k})-v(a_1)\geqslant 1$ follows from the auxiliary Lemma~\ref{lem:aux}.
\end{itemize}
Now, assume $2^{n+1}+2^n\leqslant l_0<2^{n+2}$. 
\begin{itemize}
    \item If $k\leqslant l_0$
    \begin{align*}
v(a_{p,k})-v(a_1)&\geqslant 2p-4-(l_0-k)2^{-n}+v\left(\!\binom{l_0}{k}\!\right)+s(l_0-k)\\ &\geqslant 2+s(l_0-k)-(l_0-k)2^{-n}.
\end{align*}
If $s(l_0-k)\geqslant 2$ then the last expression is positive, since $(l_0-k)2^{-n}<4$. 
 
If $s(l_0-k)=1$ then $l_0-k=2^{m}$ with $m\leqslant n+1$, so $(l_0-k)2^{-n}\leqslant 2$, and we conclude that $v(a_{p, k})-v(a_1)\geqslant 1$. If $k=l_0$ then $v(a_{p, k})-v(a_1)\geqslant 2$.
\item
If $k>l_0$, note that
\begin{align*}
    v(a_{p, k})-v(a_1)&\geqslant (k-l_0)2^{-n}+s(k)-2v(k-l_0)-2s(k-l_0)-s(l_0)+2n+2\\
    &+(p-2)(n+2)-\max_{J\in \binom{T_k}{p-2}} \sum_{j\in J}v(j)
\end{align*}
The estimate $v(a_{p, k})-v(a_1)\geqslant 1$ follows from the auxiliary Lemma~\ref{lem:aux}.
\end{itemize}
We conclude that $v(a_p)-v(a_1)\geqslant 1$ as desired.
\end{proof}

\section{Arbitrarily many Holt analytic functions sharing a root}
\label{sec:counter-examples}
Since $H_{l,\eps}(\beta)=\check{H}_{\beta, \epsilon}(l)$ is defined for~$l$ in a certain neighborhood $N_{\beta, \epsilon}(\Z_2)$of $\Zb$ in~$\Cb_2$, one may ask whether Stieltjes's conjecture could hold in a wider context, namely, whether there exists $M>0$ such that for each root $\beta$ the number of $l\in N_{\beta, \epsilon}(\Z_2)$ with $H_{l, \eps}(\beta)=0$ is bounded by~$M$.
The next theorem shows that the answer is negative.
\begin{theorem}
\label{thm:counter-example}
Suppose~$n\geqslant 2$. For each $M\in\Nb$ there exist~$M$ pairs $(l_k, \eps_k)_{k=1}^M$ and $\beta\in\bar{\Qb}_2$, $v(\beta)=2^{-n}$ with 
    $l_k\in N_{n+2}(\Z_2)$ such that $H_{l_k, \eps_k}(\beta)=0$.
    More precisely, let $H_{2^n, 0}(\beta)=0$. For each $0\leqslant j<n$ and $\eps$ there exists $l\in N_{n+2}(\Z_2)$ with $l\equiv 2^n+\lfloor j/2 \rfloor \mod 2^{n+2}$ such that
    $H_{l, \eps}(\beta)=0$.
\end{theorem}
Informally, Theorem~\ref{thm:counter-example} tells us that we can find arbitrarily long tuples of Legendre polynomials~$P_l$ sharing a common root, if we allow the degrees~$l$ to be $2$-adic numbers in~$\bar{\Qb}_2$. More precisely, it shows that for each root $\alpha$ of $P_{2^{n+1}}$, we can find at least $2n$ ``polynomials'' $P_l$ with $l\in\Ol_{\Cb_2}$ such that $P_l(\alpha)=0$.
\begin{proof}
 Let $\alpha$ be a root of $P_{2^{n+1}}$ for which $\beta=4\alpha^2/(\alpha^2-1)$ (see the definition of~$H_{l, 0}$ in~\eqref{def:Hl0}). Note (see ~\eqref{eqn:valphavbeta}) that $v(\alpha-1)=1-2^{-n-1}$ and $v(\alpha)=0$. By the Legendre recursion (Lemma~\ref{lem:legendre-recursion-rolled})
    for each $1\leqslant j<2^{n+1}$
    \begin{align*}
     (2^{n+1}+1)(2^{n+1}+2)&\ldots (2^{n+1}+j)P_{2^{n+1}+j}(\alpha)
     \\& =A_j(\alpha, 2^{n+1}) P_{2^{n+1}}(\alpha)-2^{n+1}A_{j-1}(\alpha, 2^{n+1}+1)P_{2^{n+1}-1}(\alpha)\\&=-2^{n+1}A_{j-1}(\alpha, 2^{n+1}+1)P_{2^{n+1}-1}(\alpha)
     \end{align*}
    Hence, from Lemma~\ref{lem:vP00111} below, $v(j!)+v\bigl(P_{2^{n+1}+j}(\alpha)\bigr)\geqslant n+1$, from which we conclude that for all $0\leqslant j<n$ 
    \begin{align}
    \label{ineq:bigvaluations}
       v(P_{2^{n+1}+j}(\alpha))\geqslant n+1-j+s(j)\geqslant 3.
    \end{align}
    Suppose $j$ is even. The inequality~\eqref{ineq:bigvaluations} is equivalent by  Definition~\eqref{def:Hl0} to
    \begin{align}
    \label{ineq:v(a0)even}
    v\bigl(H_{2^{n}+\lfloor j/2\rfloor, 0}(\beta)\bigr)\geqslant 3+(2^{n}+\lfloor j/2 \rfloor)2^{-n}-s(2^n+\lfloor j/2 \rfloor). 
    \end{align}
    On the other hand, if we let 
    \[F(x)=\sum_{k=0}^{\infty} \frac{\beta^k}{\binom{2k}{k}}\binom{2^n+\lfloor j/2\rfloor +2^{n+2}x}{k}^2=\sum_p a_p x^p\] then from Proposition~\ref{prop:va1} we know that $v(a_1)=3+(2^n+\lfloor j/2\rfloor)2^{-n}-s(2^n+\lfloor j/2\rfloor)$.
    Hence, the inequality~\eqref{ineq:v(a0)even} tells us that $v(a_0)\geqslant v(a_1)$.
    We conclude, by \Cref{prop:single-root}, that there exists $x\in B(0, 1)$ for which $F(x)=0$. Equivalently, there exists $l\in B(2^n+\lfloor j/2 \rfloor, 2^{-n-2})$  for which  $\check{H}_{\beta, 0}(l)=H_{l, 0}(\beta)=0$.
    The argument for the case $j$ is odd is identical to the above, once replacing the formula for $H_{2^n+\lfloor j/2\rfloor, 0}$ with the one for
    $H_{2^n+\lfloor j/2\rfloor, 1}$.
    \end{proof}

    \begin{lemma}\label{lem:vP00111}
    Let $v(\alpha-1)=1-2^{-n-1}$. Then $v\bigl(P_{2^{n+1}-1}(\alpha)\bigr)=0$
\end{lemma}
\begin{proof}
Let $\beta=4\alpha^2/(\alpha^2-1)$. Then, $v(\beta)=2^{-n}$.
Also, by definition~\eqref{def:Hl1} of~$H_{l, 1}$ we see
\begin{equation}
\label{eq:translating-valuations}
v\bigl(P_{2l+1}(\alpha)\bigr)=s(l)-l 2^{-n}+v\bigl(H_{l, 1}(\beta)\bigr).
\end{equation}
For $l=2^n-1$ we obtain
\[v\bigl(P_{2^{n+1}-1}(\alpha)\bigr)=n-1+2^{-n}+v\bigl(H_{l, 1}(\beta)\bigr).\]
Now,
\[ H_{2^n-1, 1}(\beta)=\sum_{k=0}^{2^n-1}\frac{\beta^k}{\binom{2k}{k}(2k+1)}\binom{2^n-1}{k}^2 =\sum_k h_k\ .\]
Notice that
\[v(h_k)=k2^{-n}+s(k)+2s(2^n-1-k)-2s(2^{n}-1)=k2^{-n}+s(k)+2(n-s(k))-2n=k2^{-n}-s(k).\]
The latter expression is minimal if and only if $k=2^{n}-1$, in which case it is equal to $1-2^{-n}-n$.
We get that $v\bigl(H_{l, 1}(\beta)\bigr)=1-n-2^{-n}$, and $v\bigl(P_{2^{n+1}-1}(\alpha)\bigr)=0$.
\end{proof}

\section{Auxiliary binary digits functions }\label{sec:auxiliary}

This section collects the necessary auxiliary inequalities which appeared in the text, and which concern expressions involving the sum-of-digits function. The main results needed in the above are Lemmas \ref{lem:G1-properties} and \ref{lem:aux}. 
\subsection{Binary digits function I}
Consider the following auxiliary function defined for $n, a, b\geqslant 0$.
\begin{align*}
G_1(n, a, b)&=(a+1) 2^{-n}+s(a+1+b)-2v(a+1)-2s(a+1)-s(b)+2n+2\\&=(a+1) 2^{-n}+s(a+1+b)-2s(a)-s(b)+2n
\end{align*}
\begin{lemma}[Properties of $G_1$]
The following lower bounds on~$G_1$ hold.
\label{lem:G1-properties}
        \begin{enumerate}[label=(\alph*)]
    \item\label{itm:l0-caseA}
        If $2^n\leqslant l_0<2^{n+1}$ then $G_1(n, a, l_0)\geqslant 1$ with  equality if and only if $a+1=2^n$ or $a +1 =2\cdot 2^{n}$ or $a+1=3\cdot 2^n$ or $a+1=4\cdot 2^{n}$.
        \item\label{itm:l0-caseB}
        If $2^{n+1}+2^n\leqslant l_0 < 2^{n+2}$ then $G_1(n, a, l_0)\geqslant 0$ with equality if and only if $a+1=2^n$ or $a +1 =2\cdot 2^{n}$.
    \end{enumerate}
\end{lemma}
\begin{proof}[Proof of Lemma~\ref{lem:G1-properties}]
Parts~\ref{itm:l0-caseA} and~\ref{itm:l0-caseB} both  follow from Lemma~\ref{lem:induction}
and Lemma~\ref{lem:G1-base}.
\end{proof}
\begin{lemma}[Analysis of $G_1$: induction step]
\label{lem:induction}
For $n\geqslant 1$ the inequality
    \[G_1(n, a, b)\geqslant G_1(n-1, \lfloor a/2\rfloor, \lfloor b/2\rfloor)\]
    holds, with equality if and only if $a$ is odd.
\end{lemma}
\begin{proof}
For $n\geqslant 1$, let us write $a=2a'+a_0$ and $b=2b'+b_0$ with $0\leqslant a_0, b_0\leqslant 1$.
We consider all possible cases:
\begin{enumerate}[label=(\roman*)]
\item $a_0=0$, $b_0=0$: In that case
\begin{align*}
G_1(n, a, b)&=(a'+1)2^{-n+1}+s(a'+b')+1-2s(a'+1)-s(b')+2(n-1)+(2-2^{-n})
\\&>  G_1(n-1, a', b')
\end{align*}
\item $a_0=0, b_0=1$: In that case
\begin{align*}
G_1(n, a, b)&=(a'+1)2^{-n+1}+s(a'+b'+1)-2s(a')-s(b')+2(n-1)+(1-2^{-n})\\ &>G_1(n-1, a', b')
\end{align*}
\item $a_0=1$: In that case
\begin{align*}
G_1(n, a, b)&=(2a'+2)2^{-n}+s(2a'+2b'+2+b_0)-2s(a')-2-s(b')-b_0+2n
\\&= G_1(n-1, a', b'). \qedhere
\end{align*}
\end{enumerate}
\end{proof}
\begin{lemma}[Analysis of~$G_1$: the case $n=0$]
\label{lem:G1-base}
The following  estimates on $G_1$ in the case $n=0$ hold.
\begin{enumerate}[label=(\alph*)]
\item \label{itm:g1-small-a4} If $0\leqslant a<4,$ then
\begin{enumerate}[label=(\roman*)]
    \item $G_1(0, a, 1)=a+s(a+2)-2s(a) = 1;$
    \item $G_1(0, a, 3)=a-1+s(a+4)-2s(a)\geqslant 0$ with equality if and only if $a=0$ or $a=1$.
\end{enumerate}
\item \label{itm:g1-small-a4-8} If $4\leqslant a<8,$ then
\begin{enumerate}[label=(\roman*)]
    \item $G_1(0, a, 1)\geqslant  3;$
    \item $G_1(0, a, 3)\geqslant 2.$
\end{enumerate}
    \item \label{itm:g1-large-a} If $2^A\leqslant a<2^{A+1}$ with $A\geqslant 3,$ then
\begin{enumerate}[label=(\roman*)]
    \item \label{itm:g11} $G_1(0, a, 1)\geqslant 2^A;$
    \item \label{itm:g12} $G_1(0, a, 3)\geqslant 2^A-1.$
\end{enumerate} 

\end{enumerate}
\end{lemma}
\begin{proof}
The proof of parts~\ref{itm:g1-small-a4} and~\ref{itm:g1-small-a4-8} 
is by direct calculation.
To prove part~\ref{itm:g1-large-a} write $a=2a'+a_0$ with $0\leqslant a_0\leqslant 1$.
For part~\ref{itm:g1-large-a}\ref{itm:g11} note
\begin{align*}
a+s(a+2)-2s(a)&=2a'+a_0+s(2a'+a_0+2)-2s(2a'+a_0)\\ &=2a'+a_0+s(a'+1)+a_0-2s(a')-2a_0=2a'+s(a'+1)-2s(a').
\end{align*}
The preceding expression is at least $2^{A}$ by Lemma~\ref{lem:aux2}. 
To prove part~\ref{itm:g1-large-a}\ref{itm:g12} note that due to  part~\ref{itm:g1-large-a}\ref{itm:g11} we have
\begin{align*}a+s(a+4)-2s(a)&=2a'+a_0+s(2a'+4+a_0)-2s(2a'+a_0)\\
&=a'+\bigl(a'+s(a'+2)-2s(a')\bigr)\geqslant 2^{A-1}+2^{A-1}=2^A. \qedhere 
\end{align*}

\end{proof}
\begin{lemma}
\label{lem:aux2}
Suppose $2^A\leqslant a<2^{A+1}$ with $A\geqslant 2$. Then
\[2 a+s(a+1)-2s(a)\geqslant 2^{A+1}\]
\end{lemma}
\begin{proof}
    Let $f_k(a)=2^k a+s(a+1)-2s(a)$.
    Note that $f_k(a)\geqslant f_{k+1}(\lfloor a/2 \rfloor)$: indeed,
    \begin{align*}
        f_k(2a')&=2^{k+1} a'+s(a')+1-2s(a')\geqslant f_{k+1}(a')
    \end{align*}
    and
    \begin{align*}
        f_k(2a'+1)&=2^{k+1}a'+2^k+s(a'+1)-2s(a')-2\geqslant f_{k+1}(a').
    \end{align*}
    Hence, the required estimate follows from the inequality $f_{A}(a)\geqslant 2^{A+1}$
    for $a=2, 3$. These two cases are easily verified directly.
\end{proof}
\subsection{Binary digits function II}

    Let $a, b\geqslant 0$. Let $T$ be the multiset
    \[T =[1, 1, 2, 2, 3, 3, \dots, a, a, 1, 1, 2, 2, 3, 3, \ldots, b, b]\]
        Consider the auxiliary function 
\[G_2(n, a, b, q)=q(n+2)-\max_{J\in\binom{T}{q}} \sum_{j\in J}v(j)\]

    \begin{lemma}[Analysis of $G_2$]
\label{lem:G2}
Let $1\leqslant l_0< 2^{n+2}$ and $q\geqslant 1$. 
\begin{enumerate}[label=(\alph*)]
    \item If $0\leqslant a<2^{n+2}$ then $G_2(n, a, l_0, q)\geqslant q$.
    \item If $2^{n+A}\leqslant  a<2^{n+A+1}$ where $A\geqslant 2$.
Then $G_2(n, a, l_0,q)\geqslant -2^A+2A$.
\end{enumerate}
\end{lemma}
\begin{proof}
To see the first part, observe that for all $j\in J$ we have $v(j)\leqslant n+1$.
For the second part,
we may assume that $a=2^{n+A+1}-1$ and $l_0=2^{n+2}-1$. Observe that among $1\leqslant j\leqslant a$
there are exactly $2^k$ elements of valuation $n+A-k$.
We try to make $G_2$ small by picking~$q$ elements from~$T$ of largest valuations. Observe that as long as $q$ is small enough so not all the elements in $T$ of valuation $>n+2$ are exhausted we can make $G_2$ smaller by increasing $q$.
Hence, we may assume that $q=2\sum_{k=0}^{A-3} 2^{k}=2^{A-1}-2$.
We compute
\begin{equation*}
    G_2(n, 2^{n+A+1}-1, 2^{n+2}-1, 2^{A-1}-2)=(2^{A-1}-2)(n+2)-2\sum_{k=0}^{A-3}2^k (n+A-k)=-2^A+2A.\qedhere
\end{equation*}
\end{proof}

\subsection{Binary digits function III}
Let
    \[G(n, a, b,q)\colonequals G_1(n, a, b)-2+G_2(n, a, b, q).\]
\begin{lemma} Suppose $q\geqslant 1$.
\label{lem:aux}
    \begin{enumerate}[label=(\alph*)]
    \item
        If $2^n\leqslant l_0<2^{n+1},$ then $G(n, a, l_0,q)\geqslant 0$. 
        \item
        If $2^{n+1}+2^n\leqslant l_0 < 2^{n+2},$ then $G(n, a, l_0,q)\geqslant -1$.
    \end{enumerate}
\end{lemma}

\begin{proof}
    Assume that $2^{n}\leqslant l_0<2^{n+1}$.
       According to Lemma~\ref{lem:induction} we have
    \[G_1(n, a, l_0)\geqslant G_1(0, \lfloor a2^{-n} \rfloor, 1).\]
    Hence,
    \begin{itemize}
        \item 
       If $0\leqslant a<2^{n+2}$ then according to  Lemma~\ref{lem:G1-base}\ref{itm:g1-small-a4} we have
    \[G_1(n, a, l_0)\geqslant 1.\]
     Together with Lemma~\ref{lem:G2} we obtain
    \[G(n, a, l_0,q)=G_1+G_2-2\geqslant 1+q-2=q-1.\]
\item If $2^{n+2}\leqslant a<2^{n+3}$ then by Lemma~\ref{lem:G1-base}\ref{itm:g1-small-a4-8}
and Lemma~\ref{lem:G2}
    \[G(n, a, l_0,q)=G_1+G_2-2\geqslant 3+0-2=1.\]
 \item    If  $2^{n+A}\leqslant a<2^{n+A+1}$ with $A\geqslant 3$. Then by Lemma~\ref{lem:G1-base}\ref{itm:g1-large-a} and Lemma~\ref{lem:G2}, \[G(n, a, l_0,q )=G_1+G_2-2\geqslant 2^A-2^A+2A-2=2A-2\geqslant 4.\]
    \end{itemize}
        Assume now that $2^{n+1}+2^n\leqslant l_0<2^{n+2}$.
       According to Lemma~\ref{lem:induction} we have
    \[G_1(n, a, l_0)\geqslant G_1(0, \lfloor a2^{-n} \rfloor, 3).\]
    Hence,
    \begin{itemize}
        \item 
       If $0\leqslant a<2^{n+2}$, then, according to  Lemma~\ref{lem:G1-base}\ref{itm:g1-small-a4}, we have
    \[G_1(n, a, l_0)\geqslant 0.\]
    Together with Lemma~\ref{lem:G2} we obtain
    \[G(n, a, l_0,q)=G_1+G_2-2\geqslant 0+q-2=q-2.\]
\item If $2^{n+2}\leqslant a<2^{n+3}$, then by Lemma~\ref{lem:G1-base}\ref{itm:g1-small-a4-8}
    \[G(n, a, l_0,q)=G_1+G_2-2\geqslant 2+0-2=0.\]
 \item    If  $2^{n+A}\leqslant a<2^{n+A+1}$ with $A\geqslant 3$ then \[G(n, a, l_0,q)=G_1+G_2-2\geqslant 2^A-1-2^A+2A-2=2A-3\geqslant 3. \qedhere\]
    \end{itemize}
\end{proof}

\bibliographystyle{abbrv}
\bibliography{bibliography}
\end{document}